\documentclass[11pt,a4paper]{article}
\usepackage[utf8]{inputenc}
\usepackage{amsmath}
\usepackage{amsfonts}
\usepackage{color}
\usepackage{amssymb}
\usepackage{mathrsfs}
\usepackage{esint}
\usepackage{upgreek}
\usepackage{graphicx}
\usepackage[colorinlistoftodos]{todonotes}
\usepackage[colorlinks=true]{hyperref}
\hypersetup{
	colorlinks=blue,%
	citecolor=red,%
	filecolor=black,%
	linkcolor=blue,%
	urlcolor=blue
}
\usepackage[margin=1.3in]{geometry}

\usepackage[amsmath,thmmarks,hyperref]{ntheorem}
{
	\theoremstyle{nonumberplain}
	\theoremheaderfont{\bfseries}
	\theorembodyfont{\normalfont}
	\theoremsymbol{\mbox{$\Box$}}
	\newtheorem{pf}{Proof.}
}

\numberwithin{equation}{section}
\def\R{\mathbb{R}}
\def\S{\mathbb{S}}

\def \no{\nonumber}
\def\e{\epsilon}
\def\ve{\varepsilon}
\newcommand{\ud}{\mathrm{d}}

\newtheorem{thm}{Theorem}[section]

\newtheorem{lem}{Lemma}[section]
\newtheorem{rem}{Remark}[section]

\newtheorem{prop}{\indent Proposition}[section]

\newtheorem*{thm A}{Theorem A}
\newtheorem*{thm B}{Theorem B}
\newtheorem{cor}{Corollary}[section]

\makeatletter

\newdimen\bibspace
\renewenvironment{thebibliography}[1]{%
	\section*{\refname 
		\@mkboth{\MakeUppercase\refname}{\MakeUppercase\refname}}%
	\list{\@biblabel{\@arabic\c@enumiv}}%
	{\settowidth\labelwidth{\@biblabel{#1}}%
		\leftmargin\labelwidth
		\advance\leftmargin\labelsep
		\itemsep\bibspace
		\parsep\z@skip     %
		\@openbib@code
		\usecounter{enumiv}%
		\let\p@enumiv\@empty
		\renewcommand\theenumiv{\@arabic\c@enumiv}}%
	\sloppy\clubpenalty4000\widowpenalty4000%
	\sfcode`\.\@m}
{\def\@noitemerr
	{\@latex@warning{Empty `thebibliography' environment}}%
	\endlist}

\makeatother

\makeatletter

\allowdisplaybreaks

\title{article}
\begin{document}
	\title{\bf Constrained  Moser-Trudinger-Onofri inequality and a  uniqueness criterion for the mean field equation} 
	\author{Xuezhang Chen\thanks{X. Chen: xuezhangchen@nju.edu.cn. Both authors are partially supported by NSFC (No.12271244). }~
		and Shihong Zhang\thanks{S. Zhang: dg21210019@smail.nju.edu.cn.}\\
		{\small $^{\ast}$$^{\dag}$Department of Mathematics \& IMS, Nanjing University, Nanjing 210093, P. R. China }}
		\date{}
	
		\maketitle

	
	{\noindent\small{\bf Abstract:} We establish  Moser-Trudinger-Onofri inequalities under constraint of a deviation of the second order moments from $0$, which serves as intermediate ones between Chang-Hang's inequalities under first and second order moments constraints.  A threshold for the deviation  is a uniqueness criterion for the mean field equation 
	$$-a\Delta_{\mathbb{S}^2}u+1=e^{2u} \quad \mathrm{~~on~~} \quad \mathbb{S}^2$$
when the constant $a$ is close to $\frac{1}{2}$.

		\medskip 
		
		{{\bf $\mathbf{2020}$ MSC:} 47J30, 58J37.}
		
		
		
		\section{Introduction}
		
		Recently, Moser-Trudinger-Onofri inequality under constraints has attracted a lot of attention. Chang and Hang \cite{Chang&Hang} established improved Moser-Trudinger-Onofri inequalities under higher order moments constraint. 
		
		For $m  \in \mathbb{N}$, we define
   \begin{align*}
   	\mathcal{P}_{m}=&\left\{\mathrm{all~polynomials~on~} \mathbb{R}^{3} \mathrm {~with~degree~at~most~} m\right\},\\
   	\mathring{\mathcal{P}}_{m}=&\left\{p\in \mathcal{P}_{m}\Big| \int_{\mathbb{S}^2} p \ud \mu_{\S^2}=0\right\}
\end{align*}
and
		\begin{align*}
		\mathring{\mathcal M}_1=&\left\{u\in H^1(\S^2)\Big|\int_{\S^2}x_ie^{2u}\ud V_{g_{\S^2}}=0\quad\mathrm{for}\quad i=1,2,3\right\}.
		\end{align*}

				\begin{thm A}[Chang-Hang]
		Suppose $u\in H^1(\S^2)$ and $\int_{\S^2} e^{2u}p  \ud V_{g_{\S^2}}=0$ for all  $p\in \mathring{{\mathcal{P}}}_m$, then for any $\ve>0$, there exists a positive constant $c_\ve$ such that
		\begin{align*}
			\log \fint_{\S^2} e^{2u}  \ud V_{g_{\S^2}}\leq \left(\frac{1}{4\pi N_m(\S^2)}+\ve\right)\int_{\S^2} |\nabla u|^2  \ud V_{g_{\S^2}}+2\fint_{\S^2} u  \ud V_{g_{\S^2}}+c_{\ve}.
		\end{align*}
		\end{thm A}
		
		The constant $\frac{1}{4\pi N_m(\S^2)}+\ve$ is \emph{almost sharp}. Readers are referred to \cite{Chang&Hang} for the precise definition of $N_m(\S^2)$ and its various properties. Some exact values of $N_m(\S^2)$ have been known, for instance, $N_1(\S^2)=2, N_2(\S^2)=4$ in \cite{Chang&Hang} and $N_3(\S^2)=6$ by the first author, Wei and Wu \cite[Lemma 5]{Chen-Wei-Wu}, but most  of them remain unknown. 
		
		Very recently, Chang and Gui \cite{Chang-Gui} established a sharp Moser-Trudinger type inequality, where the center of the mass allows not to be the origin.  Our first aim is to extend the above Chang-Hang type inequality by relaxing the condition $\mathring{\mathcal P}_2$.

		To continue, we introduce some notations. It sounds novel to introduce a matrix $\Lambda:=\Lambda(u)=(\Lambda_{ij})$ for each $u \in H^1(\S^2)$, where
		\begin{equation}\label{def:Lambda}
			\Lambda_{ij}=\frac{\int_{\S^2}e^{2u(x)}\left(x_i x_j-\frac{1}{3}\delta_{ij}\right)   \ud V_{g_{\S^2}}}{\int_{\S^2}  e^{2u} \ud V_{g_{\S^2}}} \qquad \mathrm{~~for~~}\quad 1 \leq i,j \leq 3		
			\end{equation}	
			can be interpreted as  a deviation of the second order moments  from $0$.
			Clearly, $\Lambda(u+c)=\Lambda(u)$ for all $c \in \R$. 	
	If we let $\|\Lambda\|=\sqrt{\sum_{i,j=1}^{3}\Lambda^2_{ij}}$, then $\|\Lambda\|$ is invariant with respect to  $O(3)$ in the sense that $\|\Lambda(u\circ A))\|=\|\Lambda(u)\|, \forall ~ A \in O(3)$,  see Proposition \ref{prop:rotation_invariant} for a simple proof.

		\begin{thm}\label{Almost sharp MTO thm}
		For every small $\delta>0$, assume $u\in \mathring{\mathcal M}_1$ with $\|\Lambda(u)\|^2\leq \frac{2}{3}-\delta$, then for any $\ve>0$ there exists a positive constant $c_{\ve,\delta}$ such that
		\begin{align}\label{ineq1:MTO_constraints}
			\log \fint_{\S^2}e^{2u}   \ud V_{g_{\S^2}}\leq \left(\frac{1}{3}+\ve\right)\fint_{\S^2} |\nabla u|^2  \ud V_{g_{\S^2}}+2\fint_{\S^2} u  \ud V_{g_{\S^2}}+c_{\ve,\delta}.
		\end{align}
		Further assume $u\in \mathring{\mathcal M}_1$ with $\|\Lambda(u)\|^2\leq \frac{1}{6}-\delta$, then for any $\ve>0$ there exists a positive constant $c_{\ve,\delta_0}$ such that
		\begin{align}\label{ineq2:MTO_constraints}
		\log \fint_{\S^2}e^{2u}   \ud V_{g_{\S^2}}\leq \left(\frac{1}{4}+\ve\right)\fint_{\S^2} |\nabla u|^2  \ud V_{g_{\S^2}}+2\fint_{\S^2} u  \ud V_{g_{\S^2}}+c_{\ve,\delta}.
		\end{align}
		\end{thm}
		
		The inequality \eqref{ineq1:MTO_constraints} serves as an intermediate case between Chang-Hang's inequalities for $m=1$ and $m=2$; the inequality \eqref{ineq2:MTO_constraints} is an extension of Chang-Hang's inequality for $m=2$ by relaxing the condition of $\mathring{\mathcal P}_2$. The numbers $
		1/3+\ve$ and $1/4+\ve$ are \emph{almost optimal}, see Proposition \ref{prop:best_cst}. Moreover,  the above thresholds $2/3$ and $1/6$ for $\|\Lambda(\cdot)\|^2$ over $\mathring{\mathcal M}_1$ are both sharp, see Proposition \ref{prop:assump_Lambda}.
		
		It is well-known that  for $a \in (0,1]$, the mean field equation
	\begin{align}\label{eq:MF}
	-a\Delta_{\S^2}u+1=e^{2u} \qquad \mathrm{~~on~~}\quad \S^2
	\end{align}
	is the Euler-Lagrange equation of the functional
	$$a\fint_{\S^2} |\nabla u|^2  \ud V_{g_{\S^2}}+2\fint_{\S^2}  u \ud V_{g_{\S^2}}-\log \fint_{\S^2} e^{2u}  \ud V_{g_{\S^2}}.$$
	It is well-known that for $a \in (0,1)$, every smooth solution of  \eqref{eq:MF} belongs to $ \mathring{\mathcal M}_1$ via the Kazdan-Warner condition.
For $a \in [1/2,1)$, the study of uniqueness theorem of \eqref{eq:MF} traced back to a conjecture raised by Chang and Yang \cite{Chang-Yang}. This conjecture was confirmed by Gui and Wei \cite{Gui-Wei} in the class of axially symmetric solutions, and had recently been completely solved by Gui and Moradifam \cite{Gui Moradifam} and  a survey therein. For $a \in [1/4,1/2)$, Shi, Sun, Tian and Wei \cite{SSTW} established that the even solution of \eqref{eq:MF} must be axially symmetry,  and every axially symmetric solution must be trivial when $a=1/3$, the latter can also be seen in \cite{Gui-Wei}. For $a \in (1/3,1/2)$, C. S. Lin \cite{Lin} had constructed nontrivial axially symmetric solutions, so was confirmed  in \cite{Gui-Hu} when $a$ is slightly strict below $1/4$. It is remarkable that when $a$ approaches to $1/3$, the uniqueness of nonzero axially symmetric solutions  is available in \cite[p.41]{SSTW}. To the best of authors' knowledge, the mean field equation for $a \in (0,1/2)$ still seems a little mysterious. See also \cite{GHMW} for a brief discussion.

		After Gui-Maradifam's work, it is natural to consider the following variational problem: For every $0<a<\frac{1}{2}$ and $c_0 \in (0,2/3)$, 	
	\begin{align}\label{prob:variational_ineq_constraint}
	\inf_{\substack{	u\in \mathring{\mathcal M}_1,\\\|\Lambda (u)\|^2\leq c_0}} \left(a\fint_{\S^2} |\nabla u|^2  \ud V_{g_{\S^2}}+2\fint_{\S^2}  u \ud V_{g_{\S^2}}-\log \fint_{\S^2} e^{2u}  \ud V_{g_{\S^2}}\right). 	\end{align}
	Clearly, the problem \eqref{prob:variational_ineq_constraint} is invariant when replacing $u$ by $u+\mathrm{constant}$. 
	Theorem \ref{Almost sharp MTO thm} enables us to  establish the existence of minimizers  whenever $a \in (1/3,1/2)$. As we shall show, for $a$ close to $1/2$,   \eqref{prob:variational_ineq_constraint}  has a `\emph{good}' Euler-Lagrange equation (differing by a constant), which is exact the mean field equation.

		We utilize a similar approach to obtain the following improved inequalities with even symmetry.
			\begin{cor}\label{Sharp MTO symmetry cor}
			Assume $u\in \mathring{\mathcal M}_1$ and $u(x)=u(-x)$. If  $\|\Lambda(u)\|^2\leq \frac{2}{3}-\delta$ for every $\delta>0$, then for any $\ve>0$,  
			\begin{align*}
			\log \fint_{\S^2}e^{2u}   \ud V_{g_{\S^2}}\leq \left(\frac{1}{4}+\ve\right)\fint_{\S^2} |\nabla u|^2  \ud V_{g_{\S^2}}+2\fint_{\S^2} u  \ud V_{g_{\S^2}}+c_{\ve,\delta}.
			\end{align*}
			Furthermore, if $\|\Lambda(u)\|^2\leq \frac{1}{6}-\delta$ for every $\delta>0$, then for any $\ve>0$, 
			\begin{align*}
			\log \fint_{\S^2}e^{2u}   \ud V_{g_{\S^2}}\leq \left(\frac{1}{6}+\ve\right)\fint_{\S^2} |\nabla u|^2  \ud V_{g_{\S^2}}+2\fint_{\S^2} u  \ud V_{g_{\S^2}}+c_{\ve,\delta}.
			\end{align*}
			
		\end{cor}
		
		In contrast to multiplicity results for the mean field equation \eqref{eq:MF} with $a \in (0,1/2)$, we wonder whether there exists a uniqueness criterion  for the mean field equation \eqref{eq:MF}. The following perturbative theorem gives an affirmative answer for $a$ close to $1/2$.
		
		\begin{thm}\label{Perturbation thm}
		There exists a positive constant $\ve_0$ such that the mean field equation
			\begin{align*}
				-a\Delta_{\S^2}u+1=e^{2u},\qquad a\in (\frac{1}{2}-\ve_0,\frac{1}{2})
			\end{align*}
			satisfying $\|\Lambda(u)\|^2\leq \frac{2}{3}-\delta$ for every $\delta>0$,
			has only a trivial solution.
		\end{thm}
		
		As an application, we prove  a sharp inequality.
		\begin{thm}\label{thm: constrained sharp ineq}
			There exists a positive constant $\ve_0$ such that for every $c_0 \in (0,\frac{2}{3})$,
			\begin{align*}
			\log \fint_{\S^2}  e^{2u} \ud V_{g_{\S^2}}\leq \left(\frac{1}{2}-\ve_0\right)\fint_{\S^2} |\nabla u|^2  \ud V_{g_{\S^2}}+2\fint_{\S^2} u  \ud V_{g_{\S^2}}
			\end{align*}
			for all  $u\in \mathring{\mathcal M}_1$ with $\|\Lambda(u)\|^2\leq c_0$.
		\end{thm}
		
		The inequality in Theorem \ref{thm: constrained sharp ineq} should be in comparison with the one in \cite[Theorem 5.1]{Chang&Hang}. Our strategy is quite different, since it is difficult to write down the Euler-Lagrange equation of \eqref{prob:variational_ineq_constraint} directly due to the inequality constraint. 
		
		The paper is organized as follows. In Section \ref{Sect2}, we prove constrained Moser-Trudinger-Onofri inequalities without/with even symmetry, and discuss the almost optimal constants in these inequalities and sharpness of the assumptions on $\|\Lambda\|$. In Section \ref{Sect3}, we first establish a Liouville type theorem of the mean field equation under constraints when $a$ is close to $1/2$, and next utilize it to establish the existence of minimizers of \eqref{prob:variational_ineq_constraint}. Besides these, we also give asymptotic behavior of solutions when $a$ goes to $1/2$ and $1/3$, respectively.

		\section{Constrained Moser-Trudinger-Onofri inequality}\label{Sect2}
		
		We first give a proof of the aforementioned property of $\Lambda$.
		\begin{prop}\label{prop:rotation_invariant}
		For $u \in H^1(\S^2)$, $\Lambda(u)$ is invariant with respect to $O(3)$, concretely, $\Lambda(u\circ A):=\Lambda(u)$ for every $A \in O(3)$.
		\end{prop}	
		\begin{pf}
		Given $u \in H^1(\S^2)$ and $A^\top=(a_{ij}) \in O(3)$, by definition \eqref{def:Lambda} of $\Lambda(u(Ax))$ we have
		\begin{align*}
		\Lambda(u(Ax))_{ij}=&\Lambda_{ij}=\frac{\int_{\S^2}e^{2u(Ax)}\left(x_i x_j-\frac{1}{3}\delta_{ij}\right)   \ud V_{g_{\S^2}}}{\int_{\S^2}  e^{2u(Ax)} \ud V_{g_{\S^2}}}\\
		=&\frac{\int_{\S^2} (a_{ik}x_k a_{jl}x_l -\frac{1}{3} \delta_{ij})e^{2u(x)} \ud V_{g_{\S^2}}}{\int_{\S^2}  e^{2u(x)} \ud V_{g_{\S^2}}}=a_{ik} \Lambda(u(x))_{kl} a_{jl},
		\end{align*}
		that is, $\Lambda(u(Ax))=A^\top \Lambda(u(x)) A$.
		Thus, it is not hard to see that
		\begin{align*}
		\|\Lambda(u(Ax))\|^2=&\mathrm{tr}\left(\Lambda(u(Ax))\Lambda(u(Ax))^\top\right)\\
		=&\mathrm{tr}\left(A^\top\Lambda(u(x))\Lambda(u(x))^\top A \right)\\
		=&\mathrm{tr}\left(\Lambda(u(x))\Lambda(u(x))^\top \right)=\|\Lambda(u(x))\|^2.
		\end{align*}
		
		\end{pf}
		
		Throughout the paper, given $u \in H^1(\S^2)$ we denote by $\overline{u}=\fint_{\S^2} u \ud V_{g_{\S^2}}$ the average of $u$ over $\S^2$. Before presenting the proof of Theorem \ref{Almost sharp MTO thm}, we state Chang and Hang's \cite[Proposition 2.4]{Chang&Hang} for readers' convenience. 
				\begin{prop}\label{prop:Chang Yang}
			Suppose $\alpha>0$, $m_i>0$, $m_i\to+\infty$, $u_i\in H^1(\S^2)$ such that $\overline{u_i}=0$, $\|\nabla u_i\|_{L^2(\S^2)}=1$ and 
			\begin{align*}
				\log\int_{\S^2} e^{2m_iu_i}  \ud V_{g_{\S^2}}\geq \alpha m_i^2.
			\end{align*}
			Also assume $u_i\rightharpoonup u$ weakly in $H^1(\S^2)$ and 
			$$|\nabla u_i|^2   \ud V_{g_{\S^2}}\to |\nabla u|^2   \ud V_{g_{\S^2}}+\sigma, \qquad
				\frac{e^{2m_iu_i}}{\int_{\S^2} e^{2m_iu_i}  \ud V_{g_{\S^2}}}\to \nu \qquad \mathrm{in ~~measure.}
			$$
			 Let 
			\begin{align*}
				\{p\in\S^2|\sigma(p)\geq 4\pi\alpha\}=\{p_1,\cdots,p_N\}.
			\end{align*}
			Then
			\begin{align*}
				\nu=\sum_{k=1}^{N}\nu_k\delta_{p_k},
			\end{align*}
			where $\nu_k \in \R_+$ and $\sum_{k=1}^{N}\nu_k=1$.
		\end{prop}

		\noindent \textbf{Proof of Theorem \ref{Almost sharp MTO thm}.}
	We argue by contradiction. Suppose there exists a sequence   $\{v_n\}\subset H^1(\S^2)$ such that $\overline{v_n}=0$, $\int_{\S^2} e^{2v_n}x_i  \ud V_{g_{\S^2}}=0$ for $i=1,2,3$ with $\|\Lambda(v_n)\|^2<2/3$, but
		\begin{align*}
			\log \fint_{\S^2} e^{2v_n}  \ud V_{g_{\S^2}}-\alpha \fint_{\S^2}  |\nabla v_n|^2 \ud V_{g_{\S^2}}\to +\infty.
		\end{align*}
		This yields  $\log \fint_{\S^2} e^{2v_n}  \ud V_{g_{\S^2}}\to +\infty$. The Moser-Trudinger-Onofri inequality
	\begin{align*}
		\log \fint_{\S^2} e^{2v_n}  \ud V_{g_{\S^2}}\leq \fint_{\S^2}  |\nabla v_n|^2 \ud V_{g_{\S^2}}
	\end{align*}
	yields $\int_{\S^2} |\nabla v_n|^2  \ud V_{g_{\S^2}}\to+\infty$. Set $m_n=\left(\int_{\S^2} |\nabla v_n|^2  \ud V_{g_{\S^2}}\right)^{1/2}$ and $u_n=\frac{v_n}{m_n}$,
	then 
	\begin{align*}
		\overline{u_n}=0\qquad \mathrm{and}\qquad \|\nabla u_n\|_{L^2(\S^2)}=1
	\end{align*}
	with
	\begin{align*}
		\log\fint_{\S^2} e^{2m_nu_n}  \ud V_{g_{\S^2}}-\frac{\alpha}{4\pi} m_n^2\to+\infty.
	\end{align*}	
		Up to a subsequence, we assume $u_n\rightharpoonup u $ weakly in $H^1(\S^2)$ and
			\begin{align*}
	|\nabla u_n|^2   \ud V_{g_{\S^2}}\to |\nabla u|^2   \ud V_{g_{\S^2}}+\sigma, \qquad	\frac{e^{2m_nu_n}}{\int_{\S^2} e^{2m_nu_n}  \ud V_{g_{\S^2}}}\to \nu
		\end{align*}
		in measure. Set
		\begin{align*}
		\{p\in\S^2|\sigma(p)\geq \alpha\}=\{p_1,\cdots,p_N\}.
		\end{align*}
	By Proposition \ref{prop:Chang Yang} we have
		\begin{align*}
		\nu=\sum_{i=1}^{N}\nu_i\delta_{p_i} \qquad \mathrm{with}\qquad \sum_{i=1}^{N} \nu_i=1.
		\end{align*}
		
		Hence, we know
		\begin{align*}
			\alpha N\leq 1,
		\end{align*}
		and
	\begin{align}\label{key formula 2}
		\sum_{i=1}^{N}\nu_i p_i=0.
	\end{align}
	
	We first assert that $N\geq N_1(\S^2)=2$. For brevity, we set $\Lambda_n:=\Lambda(v_n)$. Our discussion is divided into two cases.

	\begin{itemize}
		\item[(1)]  If there exists $\delta_0>0$ such that
		\begin{align}\label{condi 1}
		\|\Lambda_n\|^2\leq \frac{2}{3}-\delta_0,
		\end{align}
		then $N\geq 3$. If $N=2$, then \eqref{key formula 2} indicates that $p_1, p_2$ are antipodal, say $p_1=(0,0,1)$ and $p_2=-p_1$.  Meanwhile we have
	\begin{align*}
		\Lambda(v_n) \to \Lambda_{\infty}=\begin{pmatrix}
		-\frac{1}{3} &0 &0\\
		0&-\frac{1}{3}&0\\
		0&0&\frac{2}{3}
		\end{pmatrix}.
	\end{align*}
	This yields a contradiction with assumption \eqref{condi 1}. It follows that $N\geq 3$ and $\alpha\leq \frac{1}{3}$.
	\item[(2)] If there exists $\delta_0>0$ such that
	\begin{align}\label{condi 2}
	\|\Lambda(v_n)\|^2\leq \frac{1}{6}-\delta_0,
	\end{align}
	then $N\geq 4$. If $N=3$, then \eqref{key formula 2} implies that the origin and $p_i (1 \leq i \leq 3)$ lie in a plane. Without loss of generality, we assume $(p_1,p_2,p_3)\in \Pi\cap\S^2$ with $\Pi=\{x_1=0\}$, and $p_1=(0,0,1)$, $p_2=(0,w_2,z_2)$ and $p_3=(0,w_3,z_3)$, by virtue of the rotation invariant of $\|\Lambda\|$ in Proposition \ref{prop:rotation_invariant}. Similarly, we have
	\begin{align*}
		\Lambda_n\to \Lambda_{\infty}=\begin{pmatrix}
		-\frac{1}{3}&0&0\\
		0&\nu_2w_2^2+\nu_3w_3^2-\frac{1}{3}&\nu_2 w_2z_2+\nu_3 w_3z_3\\
		0&\nu_2 w_2z_2+\nu_3 w_3z_3&\nu_1+\nu_2z_2^2+\nu_3 z_3^2-\frac{1}{3}
		\end{pmatrix}.
	\end{align*}
However, 
	\begin{align}\label{Lambda infty}
		\|\Lambda_{\infty}\|^2\geq& \frac{1}{9}+\left(\nu_2w_2^2+\nu_3w_3^2-\frac{1}{3}\right)^2+\left(\nu_1+\nu_2z_2^2+\nu_3 z_3^2-\frac{1}{3}\right)^2\no\\
		\geq&\frac{1}{9}+\frac{1}{2}\left(\nu_1+\nu_2(w_2^2+z_2^2)+\nu_3(w_3^2+z_3^2)-\frac{2}{3}\right)^2\no\\
		=&\frac{1}{9}+\frac{1}{18}=\frac{1}{6},
	\end{align}
	where the last identity follows from
	\begin{align}\label{pts_on_sphere}
		w_2^2+z_2^2=w_3^2+z_3^2=1.
	\end{align}
	This is a contradiction with \eqref{condi 2}.  So we see that $N\geq 4$ and $\alpha\leq \frac{1}{4}$.
	\end{itemize}
	\hfill $\Box$
	
	\begin{rem}

When the  equalities in \eqref{Lambda infty} hold, we can determine  precise locations of concentration points $p_i (1 \leq i \leq 3)$ up to an isometry. We assume $p_1=(0,0,1)$ thanks to Proposition \ref{prop:rotation_invariant}. In this case, the matrix $\Lambda_{\infty}$ is  diagonal, whence
\begin{subequations}
\begin{align}
	\nu_2w_2^2+\nu_3w_3^2=&\frac{1}{2},\label{Rem 3.1 a0}\\
	\nu_1+\nu_2z^2_2+\nu_3z^2_3=&\frac{1}{2},\label{Rem 3.1 a01}\\
	 \nu_2w_2z_2+\nu_3w_3z_3=&0.\label{Rem 3.1 c}
\end{align}
\end{subequations}
Meanwhile, \eqref{key formula 2} gives
\begin{subequations}
\begin{align}
	\nu_2w_2+\nu_3w_3=&0,\label{Rem 3.1 a}\\
	\nu_1+\nu_2z_2+\nu_3z_3=&0.\label{Rem 3.1 b}
\end{align}
\end{subequations}
Then it follows from \eqref{Rem 3.1 a} and \eqref{Rem 3.1 c} that $\nu_2w_2(z_3-z_2)=0$. If $w_2=0$, then $w_3=0$, this directly yields $p_2=(0,0,\pm1)$ and $p_3=(0,0,\pm1)$. This is impossible as $p_i (1 \leq i \leq 3)$ are distinct. So, $w_2\not= 0$ and $z_2=z_3=z$. Using \eqref{pts_on_sphere} and $p_2\not=p_3$ we have $w_2=-w_3:=w$. By \eqref{Rem 3.1 a} we know $\nu_2=\nu_3:=\nu$. Hence, by \eqref{Rem 3.1 a0}, \eqref{Rem 3.1 a01} and \eqref{Rem 3.1 b} we have 
\begin{align*}
\nu w^2=\frac{1}{4}, \qquad \nu_1+2\nu z^2=\frac{1}{2}, \qquad  \nu_1+2\nu z=0
	\end{align*}
	and $\nu_1+2\nu=1$.
Finally, we arrive at $z=-\frac{1}{2}, w=\pm \frac{\sqrt{3}}{2}$ and $\nu_1=\nu=\frac{1}{3}$. In conclusion, we obtain
$$p_2=(0,\pm \frac{\sqrt{3}}{2},-\frac{1}{2}), \quad p_3=(0,\mp \frac{\sqrt{3}}{2},-\frac{1}{2}) \qquad \mathrm{and} \qquad \nu_1=\nu_2=\nu_3=\frac{1}{3}.$$
We \emph{\textbf{emphasize}} that each $(p_1,p_2,p_3)$ consists of the vertices of a regular triangle inscribed the unit ball.
\end{rem}

		\begin{prop}\label{prop:best_cst}
		Let $\alpha>0$ and $u\in\mathring{\mathcal M}_1$.
		\begin{enumerate}
			\item[(1)] Suppose $\|\Lambda(u)\|^2\leq \frac{2}{3}-\delta$ for every $\delta>0$,  and the following inequality with $c\in \R$
			\begin{align}\label{Almost sharp formula 1}
				\log \fint_{\S^2}  e^{2u} \ud V_{g_{\S^2}}\leq \alpha\fint_{\S^2} |\nabla u|^2  \ud V_{g_{\S^2}}+2\fint_{\S^2} u  \ud V_{g_{\S^2}}+c
			\end{align}
			holds, then $\alpha\geq \frac{1}{3}$.
			\item[(2)] Suppose $\|\Lambda(u)\|^2\leq \frac{1}{6}-\delta$ for every $\delta>0$,  and the following inequality with $c\in \R$
			\begin{align}\label{Almost sharp formula 2}
			\log \fint_{\S^2}  e^{2u} \ud V_{g_{\S^2}}\leq \alpha\fint_{\S^2} |\nabla u|^2  \ud V_{g_{\S^2}}+2\fint_{\S^2} u  \ud V_{g_{\S^2}}+c
			\end{align}
			holds, then $\alpha\geq \frac{1}{4}$.
		\end{enumerate}
	\end{prop}
	\begin{pf}
	
	Let 
		$\Phi_\ve(r)=-\log(\ve^2+r^2)$ and $\Phi_{\ve,\nu}(r)=\Phi_\ve(r)+\frac{1}{2}\log \nu$ for $\nu \in \R_+$, and $\chi$ be a smooth cut-off function in $[0,\infty)$ such that $\chi(r)=1$ in $[0,\delta]$ and $\chi(r)=0$ in $[2\delta,\infty)$. 

		For $(1)$, thanks to Proposition \ref{prop:rotation_invariant} we let $p_1=(0,0,1)$, $p_2=(0,\frac{\sqrt{3}}{2},-\frac{1}{2})$ and $p_3=(0,-\frac{\sqrt{3}}{2},-\frac{1}{2})$,  and  choose $\delta>0$ sufficiently small  such that the geodesic balls $\overline{B_{2\delta}(p_i)}\cap \overline{B_{2\delta}(p_j)}=\emptyset$ for $i\not=j$.

		We define a smooth test function by
		\begin{align*}
			u(p)=\sum_{i=1}^{3}[\chi\Phi_{\ve,\frac{1}{3}}](\stackrel{\frown}{pp_i}):=\sum_{i=1}^{3} \phi_{\ve,i}(r),
		\end{align*}
		where $r:=\stackrel{\frown}{pp_i}$ denote the geodesic distance from $p_i$ to $p$. 	
			
		First,  we have
		\begin{align*}
		\fint_{\S^2}  e^{2u} \ud V_{g_{\S^2}}=&\frac{1}{4\pi} \sum_{i=1}^{3}\int_{B_{\delta}(p_i)}e^{2\phi_{\ve,i}} \ud V_{g_{\S^2}} +O_{\delta}(1)\\
		=&\frac{1}{2}\int_0^\delta \frac{1}{(\ve^2+r^2)^2} \sin r \ud r +O_{\delta}(1)\\
		=&\frac{1}{2}\int_0^\delta \frac{r}{(\ve^2+r^2)^2}  \ud r +\int_0^\delta \frac{O(r^3)}{(\ve^2+r^2)^2}  \ud r+O_{\delta}(1)\\
		=&\frac{1}{4\ve^2}+O_{\delta}(|\log \ve|).
		\end{align*}
		
		\begin{enumerate}
		\item[(I)] We claim that $u \in \mathring{\mathcal{M}}_1$.
		
					Let
			\begin{align*}
			O(1,2)=\begin{pmatrix}
			1&0&0\\
			0&-\frac{1}{2}&\frac{\sqrt{3}}{2}\\
			0&-\frac{\sqrt{3}}{2}&-\frac{1}{2}
			\end{pmatrix} \qquad \mathrm{and}\qquad O(1,3)=\begin{pmatrix}
			1&0&0\\
			0&-\frac{1}{2}&-\frac{\sqrt{3}}{2}\\
			0&\frac{\sqrt{3}}{2}&-\frac{1}{2}
			\end{pmatrix},
			\end{align*}
			then $O(1,j) \in SO(3)$ and $O(1,j)({p}_1)={p}_j$.
			
			By the symmetry of $B_{\delta}(p_i)$ and the fact that $p_2$ and  $p_3$ are symmetric with respect to the plane $x_2=0$, we know
			\begin{align*}
			\int_{\S^2}  e^{2u}x_i \ud V_{g_{\S^2}}=0\qquad\mathrm{for}\qquad i=1,2.
			\end{align*}
			For the remaining term, we have
			\begin{align*}
			\int_{\S^2} e^{2u}x_3  \ud V_{g_{\S^2}}=&\sum_{i=1}^{3}\int_{B_{2\delta}(p_i)}e^{2u}x_3   \ud V_{g_{\S^2}}+\int_{\S^2\backslash \cup_{i=1}^3 B_{2\delta}(p_i)}x_3   \ud V_{g_{\S^2}}\\
			=&\int_{B_{2\delta}(p_1)}e^{2u}x_3   \ud V_{g_{\S^2}}+\sum_{i=2}^{3}\int_{B_{2\delta}(p_i)}e^{2u}x_3   \ud V_{g_{\S^2}}-\sum_{i=1}^{3}\int_{B_{2\delta}(p_i)}x_3   \ud V_{g_{\S^2}}\\
			=&\int_{B_{2\delta}(p_1)}\left[(e^{2u}-1)x_3 +\sum_{i=2}^{3}\left(e^{2u\circ O(1,i)}-1\right)x_3\circ O(1,i) \right ] \ud V_{g_{\S^2}}\\
			=&\int_{B_{2\delta}(p_1)} (e^{2u}-1)\left[x_3+\left(-\frac{\sqrt{3}}{2}x_2-\frac{1}{2}x_3\right)+\left(\frac{\sqrt{3}}{2}x_2-\frac{1}{2}x_3\right)\right]\ud V_{g_{\S^2}}\\
			=&0.
			\end{align*}

	\item[(II)]
		 There holds
		\begin{align}\label{step 2 formula a}
		\Lambda(u)=\begin{pmatrix}
		-\frac{1}{3}&0&0\\
		0&\frac{1}{6}&0\\
		0&0&\frac{1}{6}
		\end{pmatrix}+O_{\delta}(\ve^2|\log\ve|).
		\end{align}
		Near $p_1$, we use the  coordinates $(x_1,x_2,x_3)=(\sin r\cos\theta, \sin r\sin\theta, \cos r)$ for $r<2\delta$ and $0\leq\theta\leq 2\pi$. 
		
		Notice that
		\begin{align}\label{exp x_1^2 x_2^2}
		\frac{1}{4\pi}\int_{B_{\delta}(p_1)}e^{2u}x_1^2   \ud V_{g_{\S^2}}
		=\frac{1}{12\pi}\int_{0}^{2\pi}\int_{0}^{\delta} \frac{1}{(\ve^2+r^2)^2}\sin^3r\cos^2\theta\ud r\ud \theta 
		=O_{\delta}(|\log \ve|).
		\end{align}
		Similarly, we have
		\begin{align*}
		\int_{B_{\delta}(p_1)}e^{2u}x_2^2   \ud V_{g_{\S^2}}=O_{\delta}(|\log \ve|)
		\end{align*}	
		and
		\begin{align}\label{exp x_3^2}
		&\frac{1}{4\pi}\int_{B_{\delta}(p_1)}e^{2u}x_3^2   \ud V_{g_{\S^2}}\nonumber\\
		=&\frac{1}{6}\int_{0}^{\delta} \frac{1}{(\ve^2+r^2)^2}\cos^2r\sin r\ud r= \frac{1}{12 \ve^2}+O_{\delta}(|\log \ve|).
		\end{align}
		Notice that $x_1 \circ O(1,j)(x)=x_1$ for $j=2,3$, then
		\begin{align}\label{step 2 formula a1}
		\fint_{\S^2} e^{2u}x_1^2  \ud V_{g_{\S^2}}
		=&\frac{1}{4\pi}\sum_{i=1}^{3}\int_{B_{\delta}(p_i)}e^{2u}x_1^2   \ud V_{g_{\S^2}}+O_{\delta}(1)\nonumber\\
		=&\frac{3}{4\pi}\int_{B_{\delta}(p_1)}e^{2u}x_1^2   \ud V_{g_{\S^2}}+O_{\delta}(1)
		=O_{\delta}(|\log \ve|).
		\end{align}
		Also,
		\begin{align}\label{step 2 formula a2}
		&\fint_{\S^2} e^{2v}x_2^2  \ud V_{g_{\S^2}}\nonumber\\
		=&\frac{1}{4\pi}\sum_{i=1}^{3}\int_{B_{\delta}(p_i)}e^{2u}x_2^2   \ud V_{g_{\S^2}}+O_{\delta}(1)\nonumber\\
		=&\frac{1}{4\pi}\int_{B_{\delta}(p_1)}e^{2u}\left(x_2^2+\sum_{i=2}^{3}x^2_2\circ O(1,i)\right)   \ud V_{g_{\S^2}}+O_{\delta}(1)\nonumber\\
		=&\frac{1}{4\pi}\int_{B_{\delta}(p_1)}e^{2u}\left(x_2^2+\left(-\frac{1}{2}x_2+\frac{\sqrt{3}}{2}x_3\right)^2+\left(-\frac{1}{2}x_2-\frac{\sqrt{3}}{2}x_3\right)^2\right)   \ud V_{g_{\S^2}}+O_{\delta}(1)\nonumber\\
		=&\frac{3}{2}\cdot \frac{1}{4\pi}\left(\int_{B_{\delta}(p_1)}e^{2u}x_3^2   \ud V_{g_{\S^2}}+\int_{B_{\delta}(p_1)}e^{2u}x_2^2   \ud V_{g_{\S^2}}\right)+O_{\delta}(1)\nonumber\\
		=&\frac{1}{8 \ve^2}+O_{\delta}(|\log \ve|).
		\end{align}	
		Similarly we obtain
		\begin{align}\label{step 2 formula a3}
		\fint_{\S^2} e^{2u}x_3^2  \ud V_{g_{\S^2}}=\frac{1}{8 \ve^2}+O_{\delta}(|\log \ve|).
		\end{align}	
		Finally, it follows from the symmetry and the exact locations of $p_1,p_2,p_3$ that
		\begin{align}\label{step 2 formula a4}
			\int_{\S^2} e^{2u}x_ix_j  \ud V_{g_{\S^2}}=0\qquad\mathrm{for}\qquad i\not=j.
		\end{align}
		Hence, the claim \eqref{step 2 formula a} follows from \eqref{step 2 formula a1}, \eqref{step 2 formula a2} , \eqref{step 2 formula a3} and \eqref{step 2 formula a4}.
	\item[(III)] There hold
	\begin{align*}
		\bar{u}:=\fint_{\S^2} u  \ud V_{g_{\S^2}}=O\left(\ve^2 |\log \ve|\right)
\qquad \mathrm{and}\qquad \int_{\S^2} |\nabla u|^2  \ud V_{g_{\S^2}}=24\pi\log\frac{1}{\ve}+O_{\delta}(1).
	\end{align*}
		A direct computation yields
		\begin{align*}
			\fint_{\S^2} u  \ud V_{g_{\S^2}}=&\frac{3}{4\pi}\int_{B_{\delta}(p_1)}u\ud V_{g_{\S^2}}+O_{\delta}(1)\nonumber\\
			=&\frac{3}{2}\int_{0}^{\delta}\log\frac{1}{\ve^2+r^2}\sin r\ud r+O_{\delta}(1)\nonumber\\
			=&O_\delta\left(\ve^2 |\log \ve|\right)
		\end{align*}
		and
		\begin{align*}
		\int_{\S^2} |\nabla u|^2  \ud V_{g_{\S^2}}=& \sum_{i=1}^3 \int_{B_\delta(p_i)} |\nabla \phi_{\ve,i}|^2 \ud V_{g_{\S^2}}+O_{\delta}(1)\\
		=&6 \pi \int_0^\delta (\frac{2r}{\ve^2+r^2})^2 \sin r \ud r+O_{\delta}(1)\\
		=& 24 \pi \int_0^\delta \frac{r^3}{(\ve^2+r^2)^2} \ud r+\int_0^\delta \frac{O(r^5)}{\ve^2+r^2)^2} \ud r+O_{\delta}(1)\\
		=&24 \pi \log \frac{1}{\ve}+O_{\delta}(1).
		\end{align*}
	\end{enumerate}
	
	Therefore, fix $\delta>0$ first and choose $\ve\ll  \delta$, there hold $\|\Lambda(u)\|^2=\frac{1}{6}+O_{\delta}\left(\ve^2\log\frac{1}{\ve}\right)\leq\frac{2}{3}-\delta$ and $\int_{\S^2}  e^{2u}x_i \ud V_{g_{\S^2}}=0$ for $i=1,2,3$. Then it follows from the inequality \eqref{Almost sharp formula 1} that
		\begin{align*}
			2\log\frac{1}{\ve}\leq 6\alpha\log\frac{1}{\ve}+O\left(\ve^2 \log\frac{1}{\ve}\right).
		\end{align*}
		This directly implies $\alpha\geq \frac{1}{3}$.
		\vskip 8pt
	For $(2)$, since the strategy is very similar to $(1)$, the key ingredient is to verify $(I)$ and $(II)$.	To this end, we  select $p_1=(0,0,1)$, $p_2=(0,\frac{2\sqrt{2}}{3}, -\frac{1}{3})$, $p_3=(\sqrt{\frac{2}{3}}, -\frac{\sqrt{2}}{3},-\frac{1}{3})$ and $p_4=(-\sqrt{\frac{2}{3}}, -\frac{\sqrt{2}}{3},-\frac{1}{3})$, which are exact four vertices of a regular tetrahedron inscribed in the unit ball,  such that $\sum_{i=1}^4 p_i=0$. Let
	\begin{align*}
		O(1,2)=\begin{pmatrix}
		1&0&0\\
		0&-\frac{1}{3}&\frac{2\sqrt{2}}{3}\\
		0&-\frac{2}{\sqrt{3}}&-\frac{1}{3}
		\end{pmatrix}, \qquad
		O(1,3)=\begin{pmatrix}
		\frac{1}{\sqrt{3}}&0&\sqrt{\frac{2}{3}}\\
		\frac{2}{3}&-\frac{1}{\sqrt{3}}&-\frac{\sqrt{2}}{3}\\
		\frac{\sqrt{2}}{3}&\sqrt{\frac{2}{3}}&-\frac{1}{3}
		\end{pmatrix}
	\end{align*}
		and
		\begin{align*}
		O(1,4)=\begin{pmatrix}
		-\frac{1}{\sqrt{3}}&0&-\sqrt{\frac{2}{3}}\\
		\frac{2}{3}&\frac{1}{\sqrt{3}}&-\frac{\sqrt{2}}{3}\\
		\frac{\sqrt{2}}{3}&-\sqrt{\frac{2}{3}}&-\frac{1}{3}
		\end{pmatrix},
		\end{align*}
		then $O(1,j)(p_1)=p_j$ and $O(1,j)\in SO(3)$ for $j=2,3,4$. We define
		\begin{align}\label{fcn:test_second}
		u(p)=\sum_{i=1}^{4}\Phi_{\ve,\frac{1}{4}}(\stackrel{\frown}{pp_i}):=\sum_{i=1}^{4} \phi_{\ve,i}(r).
		\end{align}
			Near $p_1$, we choose the coordinates $(x_1,x_2,x_3)=(\sin r\cos\theta, \sin r\sin\theta, \cos r)$, where $r<2\delta$ and $0\leq\theta\leq 2\pi$. By symmetry we have
			\begin{align}\label{x near p_1}
				\int_{B_{2\delta}(p_1)}x_1   \ud V_{g_{\S^2}}=	\int_{B_{2\delta}(p_1)}x_2   \ud V_{g_{\S^2}}=0
			\end{align}
			and
			\begin{align}\label{exponential near p_1}
			\int_{B_{2\delta}(p_1)}e^{2u}x_1   \ud V_{g_{\S^2}}=	\int_{B_{2\delta}(p_1)}e^{2u}x_2   \ud V_{g_{\S^2}}=0.
			\end{align}
			Again by symmetry we have 
			$$\int_{B_{2\delta}(p_2)}e^{2u}x_1\ud V_{g_{\S^2}}=0$$
			and
			$$\int_{B_{2\delta}(p_3)}+\int_{B_{2\delta}(p_4)}e^{2u}x_1\ud V_{g_{\S^2}}=0.$$
			 Thus, we obtain $\int_{\S^2} e^{2u}x_1  \ud V_{g_{\S^2}}=0$. 
			 
			 We apply the same trick as in $(1)$ to show that
			\begin{align*}
				\int_{\S^2} e^{2u}x_2  \ud V_{g_{\S^2}}=&\sum_{i=1}^{4}\int_{B_{2\delta}(p_i)}e^{2u}x_2  \ud V_{g_{\S^2}}+\int_{\S^2\backslash \cup_{i=1}^4 B_{2\delta}(p_i)}x_2   \ud V_{g_{\S^2}}\\
				=&\int_{B_{2\delta}(p_1)}e^{2u}x_2   \ud V_{g_{\S^2}}+\sum_{i=2}^{4}\int_{B_{2\delta}(p_i)}e^{2u}x_2   \ud V_{g_{\S^2}}-\sum_{i=1}^{4}\int_{B_{2\delta}(p_i)}x_2   \ud V_{g_{\S^2}}\\
				=&\int_{B_{2\delta}(p_1)}(e^{2u}-1)x_2 +\sum_{i=2}^{4}\left(e^{2u\circ O(1,i)}-1\right)x_2\circ O(1,i)   \ud V_{g_{\S^2}}\\
				=&\int_{B_{2\delta}(p_1)}(e^{2u}-1)\left(x_2-\frac{1}{3}x_2+\frac{2\sqrt{2}}{3}x_3+	\frac{2}{3}x_1-\frac{1}{\sqrt{3}}x_2-\frac{\sqrt{2}}{3}x_3 \right.\\
			&\qquad\qquad\qquad\qquad~~\left.+	\frac{2}{3}x_1+\frac{1}{\sqrt{3}}x_2-\frac{\sqrt{2}}{3}x_3	\right)\ud V_{g_{\S^2}}\\
			=&\frac{2}{3}\int_{B_{2\delta}(p_1)}(e^{2u}-1)(2x_1+x_2)\ud V_{g_{\S^2}}\\
			=&0,
			\end{align*}
			where the last identity follows by \eqref{x near p_1} and \eqref{exponential near p_1}. Similarly, we have
			\begin{align*}
					&\int_{\S^2} e^{2u}x_3  \ud V_{g_{\S^2}}\\
					=&\int_{B_{2\delta}(p_1)}(e^{2u}-1)\left(x_3-\frac{2}{\sqrt{3}}x_2-\frac{1}{3}x_3+	\frac{\sqrt{2}}{3}x_1+\sqrt{\frac{2}{3}}x_2-\frac{1}{3}x_3+	\frac{\sqrt{2}}{3}x_1-\sqrt{\frac{2}{3}}x_2-\frac{1}{3}x_3	\right)\ud V_{g_{\S^2}}\\
					=&\frac{2}{3}\int_{B_{2\delta}(p_1)}(\sqrt{2}x_1-\sqrt{3}x_2)\ud V_{g_{\S^2}}\\
					=&0.
			\end{align*}
			So, we know $u \in \mathring{\mathcal M}_1$.
			
			Observe that
			\begin{align*}
					\int_{B_{\delta}(p_1)}x_ix_j   \ud V_{g_{\S^2}}=0 \qquad \mathrm{and}\qquad 	\int_{B_{\delta}(p_1)}e^{2u}x_jx_j   \ud V_{g_{\S^2}}=0 \qquad\mathrm{for}\qquad i\not= j.
			\end{align*}

			Similar to \eqref{exp x_1^2 x_2^2} and \eqref{exp x_3^2}, we have
			\begin{align*}
				\int_{B_{\delta}(p_1)}e^{2u}x_i^2  \ud V_{g_{\S^2}}=O_{\delta}(|\log \ve|) \qquad\mathrm{~~for~~}\qquad i=1,2
			\end{align*}
			and
			\begin{align*}
			\frac{1}{4\pi}\int_{B_{\delta}(p_1)}e^{2u}x_3^2   \ud V_{g_{\S^2}}=\frac{1}{16 \ve^2}+O_{\delta}(|\log \ve|).
			\end{align*}
			
			First, it follows from the symmetry and exact locations of $p_i$  that
			\begin{align}\label{(2) mix a1}
				\int_{\S^2} e^{2u} x_1x_i \ud V_{g_{\S^2}}=0 \qquad\mathrm{for}\qquad i=2,3.
			\end{align}
			For the last mixed term, we have
			\begin{align}\label{(2) mix a2}
					&\fint_{\S^2} e^{2u} x_2x_3 \ud V_{g_{\S^2}}\nonumber\\
					=&\frac{1}{4\pi}\sum_{i=1}^{4}\int_{B_{\delta}(p_i)} e^{2u}x_2x_3  \ud V_{g_{\S^2}}+O_{\delta}(|\log \ve|)\nonumber\\
					=&\frac{1}{4 \pi}\int_{B_{\delta}(p_1)} e^{2u}x_2x_3  \ud V_{g_{\S^2}}+\frac{1}{4\pi}\sum_{i=2}^{4}\int_{B_{\delta}(p_1)} e^{2u}(x_2 x_3)\circ O(1,i)  \ud V_{g_{\S^2}}+O_{\delta}(|\log \ve|)\nonumber\\
					=&\frac{1}{4\pi}\int_{B_{\delta}(p_i)} \left[x_2x_3 +\left(-\frac{1}{3}x_2+\frac{2\sqrt{2}}{3}x_3\right)\left(-\frac{2}{\sqrt{3}}x_2-\frac{1}{3}x_3\right)\right.\nonumber\\
					&\qquad\qquad\quad~~+\left(\frac{2}{3}x_1-\frac{1}{\sqrt{3}}x_2-\frac{\sqrt{2}}{3}x_3\right)\left(\frac{\sqrt{2}}{3}x_1+\sqrt{\frac{2}{3}}x_2-\frac{1}{3}x_3\right)\nonumber\\
					&\qquad\qquad\quad~~\left.+\left(\frac{2}{3}x_1+\frac{1}{\sqrt{3}}x_2-\frac{\sqrt{2}}{3}x_3\right)\left(\frac{\sqrt{2}}{3}x_1-\sqrt{\frac{2}{3}}x_2-\frac{1}{3}x_3\right)\right] e^{2u}  \ud V_{g_{\S^2}}+O_{\delta}(|\log \ve|)\nonumber\\
					=&O_{\delta}(|\log \ve|).
			\end{align}
			
			We next turn to the diagonal elements in $\Lambda$. Notice that
			\begin{align}\label{(2) diag a1}
				&\fint_{\S^2} e^{2u}x_1^2  \ud V_{g_{\S^2}}\nonumber\\
				=&\frac{1}{4\pi}\sum_{i=1}^{4}\int_{B_{\delta}(p_i)} e^{2u}x_1^2  \ud V_{g_{\S^2}}+O_{\delta}(1)\nonumber\\
				=&\frac{1}{4\pi}\int_{B_{\delta}(p_1)} e^{2u}\left[x_1^2+x_1^2+\left(\frac{1}{\sqrt{3}}x_1+\sqrt{\frac{2}{3}}x_3\right)^2 +\left(-\frac{1}{\sqrt{3}}x_1-\sqrt{\frac{2}{3}}x_3\right)^2\right] \ud V_{g_{\S^2}}+O_{\delta}(1)\nonumber\\
				=&\frac{1}{12\ve^2}+O_{\delta}(|\log \ve|)
			\end{align}
			and
			\begin{align}\label{(2) diag a2}
				&\fint_{\S^2} e^{2u}x_2^2  \ud V_{g_{\S^2}}\nonumber\\
				=&\frac{1}{4\pi}\sum_{i=1}^{4}\int_{B_{\delta}(p_i)} e^{2u}x_2^2  \ud V_{g_{\S^2}}+O_{\delta}(1)\nonumber\\
				=&\frac{1}{4\pi}\int_{B_{\delta}(p_1)} \left[x_2^2+\left(-\frac{1}{3}x_2+\frac{2\sqrt{2}}{3}x_3\right)^2+\left(\frac{2}{3}x_1-\frac{1}{\sqrt{3}}x_2-\frac{\sqrt{2}}{3}x_3\right)^2 \right.\nonumber\\
				&\left.\qquad\qquad\quad~~+\left(\frac{2}{3}x_1+\frac{1}{\sqrt{3}}x_2-\frac{\sqrt{2}}{3}x_3\right)^2\right]e^{2u}\ud V_{g_{\S^2}}+O_{\delta}(1)\nonumber\\
				=&\frac{1}{12\ve^2}+O_{\delta}(|\log \ve|)
			\end{align}
			and
			\begin{align}\label{(2) diag a3}
					&\fint_{\S^2} e^{2u}x_3^2  \ud V_{g_{\S^2}}\nonumber\\
					=&\frac{1}{4\pi}\sum_{i=1}^{4}\int_{B_{\delta}(p_i)} e^{2u}x_3^2  \ud V_{g_{\S^2}}+O_{\delta}(1)\nonumber\\		
				=&\frac{1}{4\pi}\int_{B_{\delta}(p_1)} \left[x_3^2+\left(-\frac{2}{\sqrt{3}}x_2-\frac{1}{3}x_3\right)^2+\left(\frac{\sqrt{2}}{3}x_1+\sqrt{\frac{2}{3}}x_2-\frac{1}{3}x_3\right)^2 \right.\nonumber\\
				&\left.\qquad \qquad\quad~~+\left(\frac{\sqrt{2}}{3}x_1-\sqrt{\frac{2}{3}}x_2-\frac{1}{3}x_3\right)^2\right]e^{2u}\ud V_{g_{\S^2}}+O_{\delta}(1)\nonumber\\
				=&\frac{1}{12\ve^2}+O_{\delta}(|\log \ve|).
			\end{align}
			
			Combining \eqref{(2) mix a1}, \eqref{(2) mix a2}, \eqref{(2) diag a1}, \eqref{(2) diag a2} and \eqref{(2) diag a3}, we obtain
			\begin{align*}
				\Lambda(u)=O\left(\ve^2|\log \ve|\right),
			\end{align*}
			This yields $\|\Lambda(u)\|^2\leq \frac{1}{6}-\delta$ for sufficiently small $\ve$. 
			
			Similar to $(1)$, it is not hard to see that
			\begin{align*}
			\overline{u}=\fint_{\S^2} u  \ud V_{g_{\S^2}}=O\left(\ve^2|\log \ve|\right)\qquad \mathrm{and}\qquad \int_{\S^2} |\nabla u|^2  \ud V_{g_{\S^2}}=32\pi|\log\ve|+O_{\delta}(1),
			\end{align*}
			
			Therefore, putting these facts together, by \eqref{Almost sharp formula 2} we obtain $\alpha\geq \frac{1}{4}$.
	\end{pf}
	
	The construction of test functions is motivated by those of Chang-Hang \cite{Chang&Hang} and Chen-Wei-Wu \cite{Chen-Wei-Wu}. Instead of the cut-paste argument used in   \cite{Chang&Hang}  or \cite{Chen-Wei-Wu},   we utilize the precise location information of $p_i$ to guarantee test functions lying in $\mathring{\mathcal M}_1$. This makes this construction more transparent. 
	
	\begin{rem}
	The test function in \eqref{fcn:test_second} can be applied to show that:
	For every small $\delta>0$, suppose $u \in $ with $\|\Lambda(u)\|^2\leq \delta$, the following inequality
	$$\log \fint_{\S^2} e^{2u}  \ud V_{g_{\S^2}} \leq \alpha\fint_{\S^2} |\nabla u|^2  \ud V_{g_{\S^2}}+2\fint_{\S^2}  u \ud V_{g_{\S^2}}+c$$
	for some $\alpha, c \in \R_+$. Then $\alpha \geq \frac{1}{4}$. In other words, even if we replace the condition $u \in \mathring{\mathcal{P}}_1$ and $\|\Lambda(u)\|\leq \frac{1}{6}-\delta$ by $u \in \mathring{\mathcal{P}}_2$, then the number $\frac{1}{4}+\ve$ can not be improved.
 	\end{rem}

	As a byproduct, we can show that the thresholds for  $\|\Lambda(\cdot)\|^2$ over $ \mathring{\mathcal M}_1$ are also sharp.

\begin{prop}\label{prop:assump_Lambda}
	For any $\alpha \in (1/3,1/2)$, there exists $\{u_n\} \subset \mathring{\mathcal M}_1$ such that 
	\begin{align}\label{Optimal Lamda 1}
		\log \fint_{\S^2} e^{2u_n}  \ud V_{g_{\S^2}}-\alpha\fint_{\S^2} |\nabla u_n|^2  \ud V_{g_{\S^2}}-2\fint_{\S^2} u_n  \ud V_{g_{\S^2}}\to+\infty
	\end{align}
	and $\|\Lambda_n\|^2\to\frac{2}{3}$, where $\Lambda_n=\Lambda(u_n)$. Similarly, for any $\alpha \in (1/4,1/3)$, there exists  $\{u_n\} \subset \mathring{\mathcal M}_1$ such that
	\begin{align}\label{Optimal Lamda 2}
	\log \fint_{\S^2} e^{2u_n}  \ud V_{g_{\S^2}}-\alpha\fint_{\S^2} |\nabla u_n|^2  \ud V_{g_{\S^2}}-2\fint_{\S^2} u_n  \ud V_{g_{\S^2}}\to+\infty
	\end{align}
	and $\|\Lambda_n\|^2\to\frac{1}{6}$.
	
\end{prop}
	\begin{pf}
	
	We adopt the same notations as in the proof of Proposition \ref{prop:best_cst}.
	
		For the first case, let $p_1=(0,0,1)$ and $p_2=(0,0,-1)$. We take
		\begin{align*}
			u_n(p)=\sum_{i=1}^{2}[\chi\Phi_{\frac{1}{n},\frac{1}{2}}](\stackrel{\frown}{pp_i}).
		\end{align*}
		A similar argument as in Proposition \ref{prop:best_cst} yields $u_n\in\mathring{\mathcal M}_1$. A direct computation gives
		\begin{align*}
		\fint_{\S^2}  e^{2u_n} \ud V_{g_{\S^2}}=\frac{n^2}{4}+O_{\delta}(\log n)
		\end{align*}
		and
				\begin{align*}
			\fint_{\S^2} u_n  \ud V_{g_{\S^2}}=O_\delta\left(n^{-2}\log n\right)\qquad \mathrm{and}\qquad \fint_{\S^2} |\nabla u|^2  \ud V_{g_{\S^2}}=4\log n+O_{\delta}(1).
		\end{align*}
		and 
		\begin{align*}
			\Lambda_n=\begin{pmatrix}
			-\frac{1}{3} &0 &0\\
			0&-\frac{1}{3}&0\\
			0&0&\frac{2}{3}
			\end{pmatrix}+O(n^{-2}\log n)~~\Longrightarrow~~ \|\Lambda_n\|^2 \to \frac{2}{3} \qquad \mathrm{~~as~~} n\to \infty.
		\end{align*}
		Then for any $\ve>0$, we combine the above estimates to derive \eqref{Optimal Lamda 1}. 
		
		For the  second case, set $p_1=(0,0,1)$, $p_2=(0,\frac{\sqrt{3}}{2},-\frac{1}{2})$ and $p_3=(0,-\frac{\sqrt{3}}{2},-\frac{1}{2})$. Similarly, it is not hard to see that
		 \begin{align*}
		u_n(p)=\sum_{i=1}^{3}[\chi\Phi_{\frac{1}{n},\frac{1}{3}}](\stackrel{\frown}{pp_i})\in \mathring{\mathcal M}_1
		\end{align*}
		and \eqref{Optimal Lamda 2} holds. Moreover, we have
		\begin{align*}
			\Lambda_n=\begin{pmatrix}
			-\frac{1}{3}&0&0\\
			0&\frac{1}{6}&0\\
			0&0&\frac{1}{6}
			\end{pmatrix}
		+O(n^{-2}\log n) ~~\Longrightarrow~~ \|\Lambda_n\|^2 \to \frac{1}{6} \qquad \mathrm{~~as~~} n\to \infty.
		\end{align*}
	\end{pf}
	
	Finally, we give a sketched proof for improved inequalities with even symmetry.
	
	\noindent \textbf{Proof of Corollary \ref{Sharp MTO symmetry cor}}: Notice that the even symmetry implies that  $u\in\mathring{\mathcal{M}}_1$, we argue by contradiction and follow the same lines in Theorem \ref{Almost sharp MTO thm}. Then it follows from Proposition \ref{prop:Chang Yang} that
	\begin{align*}
		\{p\in\S^2|\sigma(p)\geq \alpha\}=\{\pm p_1,\cdots,\pm p_m\} \quad \mathrm{~~with~~} N=2m,\quad  \alpha N \leq 1 
	\end{align*}
	and \eqref{key formula 2} with $N=2m$ holds.
	
	If $\|\Lambda\|^2\leq\frac{2}{3}-\delta_{0}$ for some $\delta_0>0$, then the same argument in Theorem \ref{Almost sharp MTO thm} gives $N\geq 3$. This together with even $N$ yields $N\geq 4$.

		 If $\|\Lambda\|^2\leq \frac{1}{6}-\delta_{0}$ for some $\delta_0>0$, then we claim that $N\geq 6$.  Suppose not, we have $N=4$ and take $p_1=(0,0,1)$, $p_2=(0, y,z)$, $p_3=(0,0,-1)$,   and $p_4=(0,-y,-z)$, where  $y\neq 0$. By \eqref{key formula 2} we have
			\begin{align*}
				\nu_2 y-\nu_4y=0\qquad \mathrm{and}\qquad \nu_1-\nu_3+\nu_2 z-\nu_4 z=0.
			\end{align*}
			This yields $\nu_1=\nu_3$ , $\nu_2=\nu_4$ and $\nu_1+\nu_2=\frac{1}{2}$. Thus, we obtain
			\begin{align*}
			\Lambda_n\to \Lambda_{\infty}=\begin{pmatrix}
			-\frac{1}{3}&0&0\\
			0&2\nu_2 y^2-\frac{1}{3}&2\nu_2 yz\\
			0&2\nu_2 yz &2\nu_1+2\nu_2 z^2-\frac{1}{3}
			\end{pmatrix}.
			\end{align*}
			However, 
			\begin{align}\label{Cor ieq a}
			\|\Lambda_{\infty}\|^2\geq& \frac{1}{9}+\left(2\nu_2 y^2-\frac{1}{3}\right)^2+\left(2\nu_1+2\nu_2 z^2-\frac{1}{3}\right)^2\nonumber\\
			\geq&\frac{1}{9}+\frac{1}{2}\left(2\nu_1+2\nu_2 (y^2+z^2)-\frac{2}{3}\right)^2\nonumber\\
			=&\frac{1}{9}+\frac{1}{18}=\frac{1}{6}.
			\end{align}
			
			This contradicts the assumption $\limsup_{n\to+\infty}\|\Lambda_n\|^2\leq \frac{1}{6}-\delta_{0}$. 
\hfill $\Box$
				
\begin{rem}
\begin{enumerate}
	\item[(1)] If the equalities in \eqref{Cor ieq a} hold, then the matrix $\Lambda_{\infty}$ is diagonal. This follows
	\begin{align*}
	2\nu_2 yz=0, \quad 2\nu_2 y^2=\frac{1}{2} \quad \mathrm{~~and~~}\quad 2\nu_1+2\nu_2 z^2=\frac{1}{2},
	\end{align*}
	so we obtain $z=0$, $y=\pm 1$ and $\nu_1=\nu_2=\nu_3=\nu_4=\frac{1}{4}$.
		\item[(2)] There exists a sequence of smooth functions  $\{u_n\}$ such that $\|\Lambda_n\|^2\to 0$, $u_n(x)=u_n(-x)$ and 
	\begin{align*}
	\lim_{n\to+\infty}\frac{\frac{1}{6}\fint_{\S^2} |\nabla u_n|^2  \ud V_{g_{\S^2}}+2\fint_{\S^2} u_n  \ud V_{g_{\S^2}}}{\log\fint_{\S^2} e^{2u_n}  \ud V_{g_{\S^2}}}=1.
	\end{align*}
	Indeed, we take
	$$u_n(p)=\sum_{i=1}^{6}[\chi\Phi_{\frac{1}{n},\frac{1}{6}}](\stackrel{\frown}{pp_i}).$$
	where $p_1=(0,0,1)$, $p_2=(0,1,0)$, $p_3=(1,0,0)$ and $p_{6-i}=-p_i$ for $1 \leq i \leq 3$. Clearly, $u_n(x)=u_n(-x)$. Other properties can be verified by similar computations in Proposition \ref{prop:best_cst}. Meanwhile, we point out that each $u_n$ is \emph{NOT} axially symmetric.
	\end{enumerate}
\end{rem}

	\section{ Uniqueness for the mean field equation under constraints}\label{Sect3}
	
	After Gui-Maradifam \cite{Gui Moradifam}, we aim to establish the uniqueness of mean field equation
	\begin{align*}
		-a\Delta_{\S^2}u+1=e^{2u}  \qquad \mathrm{~~on~~} \S^2
	\end{align*}
	with $\|\Lambda(u)\|^2\leq \frac{2}{3}-\delta$ for every $\delta>0$, where $a\in (0,\frac{1}{2})$.

Yan Yan Li  established the concentration phenomenon of the mean field equation
could occur only when $a$ is equal to $1/m$, where $m \in \mathbb{Z}_+$. The following theorem is a direction sequence of Yan Yan Li's  \cite[Theorem 0.2]{Li}.

\begin{thm B}\label{Li yanyan thm}
Suppose $a_n\to a \in (0,1)$ and $\{u_n\}$ is a sequence of smooth solutions to  
\begin{align*}
	-a_n\Delta_{\S^2}u_n+1=e^{2u_n}\qquad\mathrm{on}\qquad \S^2
\end{align*}
and $\max_{\S^2} u_n \to +\infty$. Then up to a subsequence, there exist finitely many points $p_1,\cdots, p_m$ such that
\begin{enumerate}
	\item[(1)] $u_n\to -\infty$ uniformly on $\S^2\backslash\{p_1,\cdots,p_m\}$.
	\item[(2)] There hold
	\begin{align*}
		u_n-\overline{u_n}\to 4\pi \sum_{k=1}^{m}G_{\S^2}(\cdot, p_k),
	\end{align*}
		where $G_{\S^2}(\cdot, p)$ is the Green's function of $-\Delta_{\S^2}$ with pole  $p$, and
	\begin{align*}
		e^{2u_n}\ud V_{g_{\S^2}}\rightharpoonup 4 \pi \sum_{k=1}^{m}\delta_{p_k} \qquad \mathrm{and}\qquad am=1.
	\end{align*}
\end{enumerate}
\end{thm B}	
The concentration and compactness theorem of the Gaussian curvature equation also can be found in an earlier work by Brezis and Merle \cite{Brezis-Merle}.

{\begin{thm}\label{Critical blow up thm}
		Suppose $\{u_n\}$ is a sequence of smooth solutions to
		\begin{align*}
		-a_n\Delta_{\S^2}u_n+1=e^{2u_n} \qquad\mathrm{~~with~~}\qquad a_n\to\frac{1}{3},
		\end{align*}
		then either $\{u_n\}$ are uniformly bounded or up to a rotation,
		\begin{align*}
		e^{2u_n}  \ud V_{g_{\S^2}}\rightharpoonup \frac{4 \pi}{3}(\delta_{p_1}+\delta_{p_2}+\delta_{p_3})\qquad \mathrm{~~in~~measure,}
		\end{align*}
		where $p_1=(0,0,1)$, $p_2=\left(0,\mp\frac{\sqrt{3}}{2},\mp\frac{1}{2}\right)$ and $p_3=\left(0,\pm\frac{\sqrt{3}}{2},\pm\frac{1}{2}\right)$.
	\end{thm}
	
	\begin{pf}
	Suppose the concentration phenomenon  of $\{u_n\}$ occurs, by Theorem B we have 
		\begin{align}\label{Critical case formula 1}
		e^{2u_n}\ud V_{g_{\S^2}}\rightharpoonup \frac{4\pi}{3}\sum_{k=1}^{3}\delta_{p_k}
		\end{align}
		in measure, where $p_i\in\S^2$ for $1 \leq i\leq 3$. We rotate the coordinates properly such that $p_1=(0,0,1)$. The Kazdan-Warner condition gives
		\begin{align}\label{Critical case formula 2}
		\forall~ n, \qquad \int_{\S^2}e^{2u_n}x_i   \ud V_{g_{\S^2}}=0 \qquad\mathrm{for}\qquad i=1,2,3.
		\end{align}
		Then it follows from \eqref{Critical case formula 1} and \eqref{Critical case formula 2} that
		\begin{equation}\label{p_i_linear_dep}
		\sum_{k=1}^3 p_k=0.
		\end{equation}
		This implies that $p_i (1 \leq i \leq 3)$ lie in a plane $\Pi$ passing through the origin, say $\Pi=\{x_1=0\}$. Hence, if we let $p_2=(0,w_2,z_2),p_3=(0,w_3,z_3)\in \S^2$, then the exact formulae of $p_2,p_3$ follow by \eqref{p_i_linear_dep}.
	\end{pf}
	
	The concentration phenomenon of $u_a$ for $a$ near $1/2$ is clear, see  Remark \ref{sols_near_half}. But for $a\to 1/3$, by Theorem \ref{Critical blow up thm} we know that $u_a$ is compact near $1/3$. In addition, suppose $u_{a}$ is axially symmetric and then $u_{1/3}=0$, it is not hard to show that  $u_{a}\to 0$ and $\|\Lambda(u_{a})\|^2\to 0$ as $a \to 1/3$.  Furthermore, we are able to give more delicate characterization of asymptotic behavior of $u_a$ for $a$ near $1/3$.
		
	Keep in mind that let $\{u_n\}$ be  a sequence of smooth  solutions to
	\begin{align*}
		-a_n\Delta_{\S^2}u_n+1=e^{2u_n} \qquad \mathrm{and}\qquad a_n \to \frac{1}{3},
	\end{align*}
	then we have
	
	$$\fint_{\S^2} e^{2u_n}  \ud V_{g_{\S^2}}=1,\quad u_n\in\mathring{\mathcal M}_1 \quad \Longrightarrow \quad \int_{\S^2}x u_n \ud V_{g_{\S^2}}=0.
	$$ 
	
	For brevity, we assume $u_n=u_n(x_3)$. Then
	\begin{align*}
		\Lambda(u_n)=\mathrm{diag}\{
		\Lambda_{11},\Lambda_{22},\Lambda_{33}\}
		 \qquad \mathrm{and}\qquad \Lambda_{11}=\Lambda_{22}=-\frac{1}{2}\Lambda_{33}.
	\end{align*}
	We write $\beta_n:=\Lambda_{33}(u_n)=\fint_{\S^2} e^{2u_n}\left(x_3^2-\frac{1}{3}\right)\ud V_{g_{\S^2}}$
	 and define
	\begin{align*}
		\widehat{u}_n:=&u_n-\overline{u_n}-\frac{45}{4}\fint_{\S^2}u_n\left(x_3^2-\frac{1}{3}\right)\ud V_{g_{\S^2}}\left(x_3^2-\frac{1}{3}\right)\\
		=&u_n-\overline{u_n}-\frac{15}{8a_n}\beta_n \left(x_3^2-\frac{1}{3}\right).
	\end{align*}
	\begin{thm}
	Suppose $a_n\to 1/3$ and $u_n=u_n(x_3)$ is a smooth solution of 
	\begin{align}\label{near 1/3 equ}
		-a_n\Delta_{\S^2}u_n+1=e^{2u_n},
	\end{align}
	then $u_n\to 0$ in $C^{\infty}(\S^2)$,
	\begin{align*}
		u_n(x)=\|u_n\|_{L^{\infty}(\S^2)}\left[\frac{3}{2}\left(x_3^2-\frac{1}{3}\right)+o_n(1)\right]
	\end{align*}
	and
	\begin{align*}
	\beta_n=\frac{4}{15}\|u_n\|_{L^{\infty}(\S^2)}(1+o_n(1)).
	\end{align*}
	Moreover, there holds
\begin{align*}
	\|\widehat{u}_n\|_{L^2(\S^2)}\leq C\|u_n\|^2_{L^{\infty}(\S^2)}.
\end{align*}
\end{thm}
\begin{pf}
Once the concentration phenomenon of  $\{u_n\}$ could occur, there only admits two possible concentration points due to axial symmetry, however this is impossible by virtue of Theorem  \ref{Critical blow up thm}. Hence, there exists a positive constant $C$ independent of $n$ such that 
	\begin{align*}
		\|u_n\|_{L^{\infty}(\S^2)}\leq C.
	\end{align*} 
	The standard elliptic estimates give
	\begin{align*}
		\|u_n\|_{H^2(\S^2)}\leq C\left(\|u_n\|_{L^2(\S^2)}+\|e^{2u_n}-1\|_{L^2(\S^2)}\right)
		\leq C\|u_n\|_{L^{2}(\S^2)}.
	\end{align*}
	By Sobolev embedding theorem we conclude that for any $k\in \mathbb{N}$,
	\begin{align}\label{elliptic estimate}
		\|u_n\|_{C^{k}(\S^2)}\leq C(k)\|u_n\|_{L^{\infty}(\S^2)}.
	\end{align}
	Taking limit in  \eqref{near 1/3 equ} we obtain $u_n\to u_{\infty}$ in $C^3 (\S^2)$, where $u_{\infty}=u_{\infty}(x_3)$ satisfies
	\begin{align*}
		-\frac{1}{3}\Delta_{\S^2}u_{\infty}+1=e^{2u_{\infty}} \qquad \Longrightarrow\qquad u_{\infty}=0.
	\end{align*}
	
		Let $v_n=\frac{u_n}{\|u_n\|_{L^{\infty}(\S^2)}}$, then $v_n$ satisfies
	\begin{align}\label{rewrite near 1/3 equ}
		-\frac{1}{3}\Delta_{\S^2}v_n+\left(\frac{1}{3}-a_n\right)\Delta_{\S^2}v_n=\frac{e^{2u_n}-1}{2u_n}2v_n.
	\end{align}
	The elliptic estimates together with \eqref{elliptic estimate} imply that
	\begin{align*}
	\|v_n\|_{C^k(\S^2)}\leq C(k).
	\end{align*}
	Up to a subsequence, taking the limit in \eqref{rewrite near 1/3 equ} we have $v_n\to v_{\infty}$ in $C^3(\S^2)$ and 
	\begin{align*}
		-\frac{1}{3}\Delta_{\S^2}v_{\infty}=2v_{\infty} \qquad \mathrm{with}\qquad v_{\infty}=v_{\infty}(x_3).
	\end{align*}
	Noticing that $\|v_{\infty}\|_{L^{\infty}(\S^2)}=1$, we conclude that $v_{\infty}(x_3)=\frac{3}{2}\left(x_3^2-\frac{1}{3}\right)$.
	
	Write
	\begin{align*}
			u_n(x_3)=\|u_n\|_{L^{\infty}(\S^2)}\left[\frac{3}{2}\left(x_3^2-\frac{1}{3}\right)+o_n(1)\right],
	\end{align*}
then
	\begin{align*}
		\beta_n=&6a_n\fint_{\S^2}u_n\left(x_3^2-\frac{1}{3}\right)\ud V_{g_{\S^2}}\\
		=&9a_n\|u_n\|_{L^{\infty}(\S^2)}(1+o_n(1))\fint_{\S^2}\left(x_3^2-\frac{1}{3}\right)^2\ud V_{g_{\S^2}}\\
		=&\frac{4}{15}\|u_n\|_{L^{\infty}(\S^2)}(1+o_n(1)).
	\end{align*}
	By definition of $\widehat{u}_n$ we have
	\begin{align}\label{appendix formula 1}
		12\int_{\S^2}\widehat{u}_n^2   \ud V_{g_{\S^2}}\leq \int_{\S^2} \left(-\Delta\widehat{u}_n\right) \widehat{u}_n  \ud V_{g_{\S^2}}=\frac{1}{a_n}\int_{\S^2}\left(e^{2u_n}-1\right)\widehat{u}_n \ud V_{g_{\S^2}}.
	\end{align}
	Using 
	$$u_n=\widehat{u}_n+\overline{u_n}+\frac{15}{8a_n}\beta_n \left(x_3^2-\frac{1}{3}\right),$$
	we expand
	\begin{align}\label{appendix formula 2}
		&e^{2u_n}-1\nonumber\\
		=&2\widehat{u}_n+2\overline{u_n}+\frac{15}{4a_n}\beta_n\left(x_3^2-\frac{1}{3}\right)+2\widehat{u}_n^2+2\overline{u_n}^2
		+\frac{225}{32a^2_n}\beta_n^2\left(x_3^2-\frac{1}{3}\right)^2\nonumber\\
		&+4\widehat{u}_n\overline{u_n}+\frac{15}{2a_n}\beta_n\widehat{u}_n\left(x_3^2-\frac{1}{3}\right)+\frac{15}{2a_n}\beta_n\overline{u_n}\left(x_3^2-\frac{1}{3}\right)+O\left(\|u_n\|^3_{L^{\infty}(\S^2)}\right).
	\end{align}
	Notice that \begin{align*}
		\int_{\S^2}\widehat{u}_n   \ud V_{g_{\S^2}}=	\int_{\S^2}\widehat{u}_n\left(x_3^2-\frac{1}{3}\right)  \ud V_{g_{\S^2}}=0
	\end{align*}
	and \begin{align}\label{appendix formula 3}
		\overline{u_n}=o_n(1)\|u_n\|_{L^{\infty}(\S^2)}\qquad|\widehat{u}_n|\leq C\|u_n\|_{L^{\infty}(\S^2)}.
	\end{align}
	
	We combine \eqref{appendix formula 1}, \eqref{appendix formula 2} and \eqref{appendix formula 3} to show
	\begin{align*}
		12\int_{\S^2}\widehat{u}_n^2   \ud V_{g_{\S^2}}
		\leq \frac{2}{a_n}\int_{\S^2}\widehat{u}_n^2 \ud V_{g_{\S^2}}+C\int_{\S^2}(\|u_n\|^2_{L^{\infty}(\S^2)} |\widehat{u}_n|+\|u_n\|_{L^{\infty}(\S^2)} |\widehat{u}_n|^2)\ud V_{g_{\S^2}}.
	\end{align*}
	By H\"older inequality we obtain
	\begin{align*}
		\|\widehat{u}_n\|_{L^2(\S^2)}\leq C\|u_n\|^2_{L^{\infty}(\S^2)}.
	\end{align*}
\end{pf}

As a whole, we outline the proofs of Theorems \ref{Perturbation thm} and  \ref{thm: constrained sharp ineq} consisting of three steps.
	
	 \emph{Step 1.}
		For every $\delta>0$ and $a \in (1/3,1/2)$  the full set of smooth solutions to
		\begin{align*}
		 -a\Delta_{\S^2}u_a+1=e^{2u_a}\quad\mathrm{~~and~~}\quad \|\Lambda(u_a)\|^2 \leq \frac{2}{3}-\delta
		\end{align*}
		is compact in $C^{\infty}(\S^2)$ topology.
			
		\emph{Step 2.}  To establish a Liouville-type theorem for the mean field equation under constraints, that is, Theorem \ref{Perturbation thm}.

\emph{Step 3.}
	There exists a positive constant $\ve_0$ such that for every $a\in (1/2-\ve_0,1/2)$ and $\delta>0$, the minimizer of
\begin{align*}
\min_{\substack{u\in \mathring{\mathcal M}_1\\\|\Lambda (u)\|^2\leq \frac{2}{3}-\delta}} \left(a\fint_{\S^2} |\nabla u|^2  \ud V_{g_{\S^2}}+2\fint_{\S^2}  u \ud V_{g_{\S^2}}-\log \fint_{\S^2} e^{2u}  \ud V_{g_{\S^2}}\right)
\end{align*}
satisfies
\begin{align*}
-a\Delta_{\S^2}u+1=e^{2u}.
\end{align*}
This follows Theorem \ref{thm: constrained sharp ineq}.

	\begin{thm}[Compactness]\label{Compactness Theorem}
		For every small $\delta>0$ and $\ve_0$ the set
			\begin{align*}
		\bigcup_{a \in [1/3+\ve_0,1/2)}\left\{u_{a}\Big| -a\Delta_{\S^2}u_a+1=e^{2u_a},~~\|\Lambda(u_a)\|^2\leq \frac{2}{3}-\delta\right\} 
		\end{align*}
		is compact in $C^{\infty}(\S^2)$ topology.
	\end{thm}
	\begin{pf}
	 If not, there exist some small $\delta_0>0$, a sequence $\{u_n\} $  and $a_n\to a$ such that 
	\begin{align*}
		\max_{\S^2}u_n=+\infty\qquad \mathrm{and}\qquad \|\Lambda(u_n)\|^2\leq \frac{2}{3}-\delta_{0}.
	\end{align*}
	It follows from Theorem B that there exist distinct points $p_1,\cdots, p_m$ such that 
	\begin{align*}
	e^{2u_n}\ud V_{g_{\S^2}}\rightharpoonup 4\pi a\sum_{k=1}^{m}\delta_{p_k} \qquad \mathrm{and}\qquad am=1.
	\end{align*}

\begin{enumerate}
	\item[(1)] $m\geq 2$. 
	
	If $m=1$, then there exists $p_1$ such that 
		\begin{align}\label{compact a_1}
	e^{2u_n}\ud V_{g_{\S^2}}\rightharpoonup 4\pi \delta_{p_1}.
	\end{align}
	Without loss of generality, we assume  $p_1=(0,0,1)$. The
	  Kazdan-Warner condition implies
	\begin{align*}
	\forall ~ n, \quad \int_{\S^2}  e^{2u_n}x \ud V_{g_{\S^2}}=0.
	\end{align*}
	However, this contradicts 
	\begin{align*}
		\int_{\S^2}  e^{2u_n}x_3\ud V_{g_{\S^2}}\to 4\pi \qquad \mathrm{as}\qquad n\to+\infty
	\end{align*}
	by virtue of \eqref{compact a_1}.
	\item[(2)] $m\geq 3$. 
	
	If $m=2$, then there exist $p_1$ and $p_2$ such that 
	\begin{align*}
	e^{2u_n}\ud V_{g_{\S^2}}\rightharpoonup 2\pi(\delta_{p_1}+\delta_{p_2}).
	\end{align*}
	Similarly, we  assume $p_1=(0,0,1)$. Since $u \in  \mathring{\mathcal M}_1$, we have $p_2=-p_1=(0,0,-1)$. Then
\begin{align*}
\Lambda(u_n)\to \Lambda_{\infty}= \begin{pmatrix}
-\frac{1}{3} &0 &0\\
0&-\frac{1}{3}&0\\
0&0&\frac{2}{3}
\end{pmatrix}.
\end{align*}
This clearly yields a contradiction with $\|\Lambda(u_n)\|^2\leq \frac{2}{3}-\delta_{0}$ for all $n$.
\end{enumerate}

Hence, we conclude that $m\geq 3$ together with $am=1$, whence $a\leq \frac{1}{3}$. This contradicts the assumption  $a>\frac{1}{3}$.
	\end{pf}

The above approach can be further improved to show 
\begin{rem}
		For every small $\delta>0$ and $\ve_0$ the set
			\begin{align*}
		\bigcup_{a \in [1/4+\ve_0,1/2)}\left\{u_{a}\Big| -a\Delta_{\S^2}u_a+1=e^{2u_a},~~\|\Lambda(u_a)\|^2\leq \frac{1}{6}-\delta\right\} 
		\end{align*}
		is compact in $C^{\infty}(\S^2)$ topology.
	\end{rem}

	We start with the proof of Theorem \ref{Perturbation thm}.
	
\noindent\textbf{Proof of Theorem \ref{Perturbation thm}}: We argue by contradiction. Suppose there exist $a_n\nearrow\frac{1}{2}$ and $u_n\not\equiv 0$  such that
\begin{align*}
-a_n\Delta_{\S^2}u_{n}+1=e^{2u_{n}}.
\end{align*}
	By Theorem \ref{Compactness Theorem}, we know $u_n\to u_{\infty}$ in $C^{3}(\S^2)$ and $u_{\infty}$ satisfies
	\begin{align*}
		-\frac{1}{2}\Delta_{\S^2}u_{\infty}+1=e^{2u_{\infty}}.
	\end{align*}
	By Gui-Moradifam \cite{Gui Moradifam}, we have $u_{\infty}= 0$. Set $\widehat{u}_n=u_n-\overline{u_n}\neq 0$, then
	\begin{align*}
	\int_{\S^2}\widehat{u}_n\ud V_{g_{\S^2}}=0\qquad\mathrm{and}\qquad u_n \in \mathring{\mathcal{M}}_1 ~~\Longrightarrow~~ \int_{\S^n}x\widehat{u}_n\ud V_{g_{\S^2}}=0.
	\end{align*}
	This follows that
	\begin{align*}
	6\int_{\S^2}\widehat{u}_n^2\ud V_{g_{\S^2}}\leq& \int_{\S^2}\widehat{u}_n \left(-\Delta_{\S^2}\widehat{u}_n\right)\ud V_{g_{\S^2}}=\frac{1}{a_n}e^{2\overline{u_n}}\int_{\S^n}\left(e^{2\widehat{u}_n}-1\right)\widehat{u}_n\ud V_{g_{\S^2}}\nonumber\\
	=&4(1+o_n(1))\fint_{\S^n}\widehat{u}_n^2\ud V_{g_{\S^2}},
	\end{align*}
	where the last equality follows from the fact that  $\overline{u_n}\to 0$ , $\widehat{u}_n\to 0$ uniformly on $\S^2$ and
	\begin{align*}
	e^{2\widehat{u}_n}-1=2(1+o_n(1))\widehat{u}_n.
	\end{align*} 
Hence, for all sufficiently large $n$, $\widehat{u}_n=0$ and thus $u_n=0$. However, it contradicts our assumption $u_{n}\not\equiv 0$.	
\hfill $\Box$
	
	\begin{rem}\label{sols_near_half}
 It is remarkable that C. S.	 Lin \cite{Lin} had constructed nontrivial axially symmetric solutions, denoted  by $u_a$, to 
	\begin{align*}
		-a\Delta_{\S^2}u+1=e^{2u}\qquad \mathrm{~~on~~} \quad \S^2
	\end{align*}
	for every $a\in (1/3,1/2)$.
	We would like to point out that if $\|\Lambda(u_a)\|^2\to 2/3$ as $a\to1/2$, then  $u_a$ must blow up. Otherwise, if there exists a positive constant $C$ such that $|u_a|\leq C$ when $a$ close to $1/2$, then the above proof implies  $u_a= 0$ provided that $|a-1/2|\ll  1$. This is a contradiction with Lin's construction. Hence, a direct consequence is that $\max_{\S^2}u_a\to+\infty$ as $a\to1/2$. A similar argument in case $(2)$ in the proof of Theorem \ref{Compactness Theorem} yields $\|\Lambda(u_a)\|^2\to2/3$. So, the constraint  $\|\Lambda(u)\|^2\leq 2/3-\delta$ for the uniqueness is sharp. 
	\end{rem}

	For $a \in (0,1]$, we define
	\begin{align*}
		S_a[u]=a\fint_{\S^2} |\nabla u|^2  \ud V_{g_{\S^2}}+2\fint_{\S^2}  u \ud V_{g_{\S^2}} \quad \mathrm{and}\quad S[u]:=S_1[u], ~~\forall~ u \in H^1(\S^2),
	\end{align*}
	and for every $\delta>0$,
	\begin{align}\label{prob:constrained_minimizing}
		J_a:=\inf_{\substack{u\in \mathring{\mathcal M}_1\\ \|\Lambda (u)\|^2\leq\frac{2}{3}-\delta}}\left(S_a[u] -\log \fint_{\S^2} e^{2u}  \ud V_{g_{\S^2}}.\right).
	\end{align}
	Clearly, it follows from Theorem \ref{Almost sharp MTO thm} that the above constrained variational functional has a lower bound for $a \in (1/3, 1]$.
	
	\begin{lem}\label{Vital lem}
	Suppose $u\in\mathring{\mathcal M}_1$, $\fint_{\S^2} e^{2u}  \ud V_{g_{\S^2}}=1$, $\|u\|_{H^1(\S^2)}\leq K_1$ and $\|\Lambda(u)\|^2\geq K_2>0$, then there exists a positive constant $C_0=C_0(K_1,K_2)$ such that
	\begin{align*}
		\frac{1}{2}\fint_{\S^2} |\nabla u|^2  \ud V_{g_{\S^2}}+2\fint_{\S^2} u  \ud V_{g_{\S^2}} \geq C_0.
	\end{align*}
	Moreover, there exists $\ve_0=\ve_0(K_1, K_2)$ such that 
	\begin{align*}
	 \left(\frac{1}{2}-\ve_0\right)\fint_{\S^2} |\nabla u|^2  \ud V_{g_{\S^2}}+2\fint_{\S^2} u  \ud V_{g_{\S^2}} \geq \frac{C_0}{2}.
	\end{align*}
\end{lem}
	\begin{pf}
		Suppose not, there exists a sequence $\{u_n\}$ such that $\|u_n\|_{H^1(\S^2)}\leq K_1$,  $\|\Lambda(u_n)\|^2\geq K_2$ and 
		\begin{align*}
		0\leq \frac{1}{2}\fint_{\S^2} |\nabla u_n|^2  \ud V_{g_{\S^2}}+2\fint_{\S^2} u_n  \ud V_{g_{\S^2}}\to 0.
		\end{align*}
		Up to a subsequence, $u_n\rightharpoonup u$ weakly in $H^1(\S^2)$ and  
		$$\fint_{\S^2} e^{2u}  \ud V_{g_{\S^2}}=1, \quad \|\Lambda(u)\|^2\geq K_2>0$$ 
		in virtue of the classical Moser-Trudinger inequality and Vitali convergence theorem.  Then 
		\begin{align*}
			 \frac{1}{2}\fint_{\S^2} |\nabla u|^2  \ud V_{g_{\S^2}}+2\fint_{\S^2} u  \ud V_{g_{\S^2}}\leq \liminf_{n\to+\infty} \left(\frac{1}{2}\fint_{\S^2} |\nabla u_n|^2  \ud V_{g_{\S^2}}+2\fint_{\S^2} u_n  \ud V_{g_{\S^2}}\right),
		\end{align*}
	gives
		\begin{align*}
		\frac{1}{2}\fint_{\S^2} |\nabla u|^2  \ud V_{g_{\S^2}}+2\fint_{\S^2} u  \ud V_{g_{\S^2}}\leq 0.
		\end{align*}
		It follows from Gui and Moradifam \cite{Gui Moradifam} that  $u=0$. This together the fact $\|\Lambda(0)\|^2=0$ gives us a contradiction.

		The second assertion directly follows from the first one.
	\end{pf}

\noindent\textbf{Proof of Theorem \ref{thm: constrained sharp ineq}}: First we fix $a\in \left(\frac{5}{12},\frac{1}{2}\right]$ for simplicity. Let $\{u_n\}$ be a minimizing sequence for \eqref{prob:constrained_minimizing} such that $\fint_{\S^2} e^{2u_n}  \ud V_{g_{\S^2}}=1$, $u_n\in\mathring{\mathcal M}_1
	$, $\|\Lambda(u_n)\|^2\leq \frac{2}{3}-\delta$ for every small $\delta>0$, and
	\begin{align*}
	S_a[u_n]\to J_a\qquad \mathrm{and}\qquad J_a\leq S_a[0]=0.
	\end{align*}
	
		 By Theorem \ref{Almost sharp MTO thm} we know that for any $\ve>0$, there exist $c_{\ve,\delta}$ such that
	\begin{align*}
			0\leq \left(\frac{1}{3}+\ve\right)\fint_{\S^2} |\nabla u_n|^2  \ud V_{g_{\S^2}}+2\fint_{\S^2} u_n  \ud V_{g_{\S^2}}+c_{\ve,\delta}.
	\end{align*}
	Taking $\ve$ small enough, with a positive constant $C$ independent of $a$ we have
	\begin{align}\label{gradient energy}
		\fint_{\S^2} |\nabla u_n|^2  \ud V_{g_{\S^2}}\leq C\left(S_a[u_n]+c_{\ve,\delta}\right)\leq C.
	\end{align}
	On the one hand, by Jensen inequality we have
	\begin{align*}
		2\overline{u_n} \leq \log \fint_{\S^2} e^{2u_n}  \ud V_{g_{\S^2}}=0.
	\end{align*}
On the other hand, the Moser-Trudinger-Onofri inequality gives
\begin{align*}
	2\overline{u_n}\geq -\fint_{\S^2} |\nabla u_n|^2  \ud V_{g_{\S^2}}\geq -C.
\end{align*}
These together give the bound of $|\overline{u_n}|$.

By Raleigh inequality we have
\begin{align*}
	2\int_{\S^2}\left(u_n-\overline{u_n}\right)^2\ud V_{g_{\S^2}}\leq \int_{\S^2} |\nabla u_n|^2  \ud V_{g_{\S^2}}\leq C,
\end{align*}
whence,
\begin{align}\label{L 2 energy}
	\int_{\S^2}u_n^2   \ud V_{g_{\S^2}}\leq C.
\end{align}
Thus, $\{u_n\}$ is bounded in $H^1(\S^2)$. Up to a subsequence, we assume $u_n\rightharpoonup u_a$ weakly in $H^1(\S^2)$ and then for all $p>0$, $e^{pu_n}  \to  e^{pu_a}$ in $L^1(\S^2)$ by Moser-Trudinger inequality and Vitali convergence theorem.
Hence, we obtain  
\begin{align*}
	\fint_{\S^2} e^{2u_a}  \ud V_{g_{\S^2}}=1, \qquad u_a\in \mathring{\mathcal M}_1 \qquad \mathrm{and}\qquad \|\Lambda(u_a)\|^2\leq \frac{2}{3}-\delta
\end{align*}
and also
\begin{align*}
J_a\leq S_a[u_a]\leq \liminf_{n \to \infty}S_{a}[u_n]=J_a.
\end{align*}
This indicates that $u_a$ is a minimizer for \eqref{prob:constrained_minimizing}. Moreover,
by \eqref{gradient energy} and \eqref{L 2 energy}, there exists a positive constant $C$ independent of $a$ such that 
$$\|u_a\|_{H^1(\S^2)}\leq C.$$

For convenience we let
$$\mathcal{S}:=\mathring{\mathcal{M}}_1 \cap \left\{u \in H^1(\S^2)\big|~~ \fint_{\S^2} e^{2u}  \ud V_{g_{\S^2}}=1 \right\}.$$

We next claim that 
\begin{align*}
	\|\Lambda(u_a)\|^2\to0 \qquad\mathrm{~~as~~}\quad a\nearrow\frac{1}{2}.
\end{align*}

Suppose not, there exists a sequence $a_k\nearrow\frac{1}{2}$ such that 
\begin{align*}
u_{a_k} \in \mathcal{S},\qquad	\|u_{a_k}\|_{H^1(\S^2)}\leq C \qquad \mathrm{and}\qquad \|\Lambda(u_{a_k})\|^2\geq K>0.
\end{align*}
It follows from Lemma \ref{Vital lem} that for all sufficiently large $k$,
\begin{align*}
	0<C(\delta,K)\leq a_k\fint_{\S^2} |\nabla u_{a_k}|^2  \ud V_{g_{\S^2}}+2\fint_{\S^2}  u_{a_k} \ud V_{g_{\S^2}},
\end{align*}
which contradicts
\begin{align*}
	0<C(\delta_{0},K)\leq S_{a_k}[u_{a_k}]\leq S_{a_k}[0]=0.
\end{align*} 
Moreover, we can find a small $\ve_0>0$ such that  for all $a\in \left(\frac{1}{2}-\ve_0,\frac{1}{2}\right)$,
\begin{align}\label{Lamda u_a value}
		u_a\in \mathcal{S} \qquad \mathrm{and}\qquad \|\Lambda(u_a)\|^2\leq \frac{1}{2}\left(\frac{2}{3}-\delta\right).
\end{align}

Fix $a\in \left(\frac{1}{2}-\ve_0,\frac{1}{2}\right)$ and  $\varphi\in C^{\infty}(\S^2)$. We define
\begin{align*}
	e^{2u(t)}:=&e^{2(u_a+t\varphi)}-3\left(\fint_{\S^2}e^{2(u_a+t\varphi)}x_i  \ud V_{g_{\S^2}}\right)x_i\\
	=&e^{2u_a}\left(1+O(t)\right)-3\left(\fint_{\S^2}e^{2u_a}\left(1+O(t)\right)x_i  \ud V_{g_{\S^2}}\right)x_i\\
	=&e^{2u_a}\left(1+O(t)\right).
\end{align*} 
Hence for all $|t|\ll  1$, $u(t)$ is well defined and $u(t)\in \mathcal{S}$. 

Observe that
\begin{align*}
	\frac{\ud u(t)}{\ud t}\Big|_{t=0}=\varphi-3\left(\fint_{\S^2}\varphi e^{2u_a}x_i  \ud V_{g_{\S^2}}\right)e^{-2u_a}x_i.
\end{align*}
Combining \eqref{Lamda u_a value}, we know for  $|t|\ll  1$, $u(t)$ is a smooth path in $\mathring{\mathcal M}_1$ with  $\|\Lambda(u(t))\|^2\leq \frac{2}{3}-\delta$.

 Since $u_a$ is a minimizer, we obtain
\begin{align*}
0=&\frac{\ud}{\ud t}\Big|_{t=0}\left[S_a[u(t)] -\log \fint_{\S^2} e^{2u(t)}  \ud V_{g_{\S^2}}\right]\\
=&2\fint_{\S^2} \left(-a \Delta_{\S^2} u_a+1\right)\frac{\ud u(t)}{\ud t}\Big|_{t=0}\ud V_{g_{\S^2}}-2\fint_{\S^2}e^{2u_a}\frac{\ud u(t)}{\ud t}\Big|_{t=0}\ud V_{g_{\S^2}}\\
=&2\fint_{\S^2} \left(-a \Delta_{\S^2} u_a+1\right)\varphi\ud V_{g_{\S^2}}-2\fint_{\S^2} e^{2u_a} \varphi\ud V_{g_{\S^2}}\\
&-6[\fint_{\S^2} \left(-a\Delta_{\S^2} u_a+1\right)e^{-2u_a}x_i\ud V_{g_{\S^2}}]
\fint_{\S^2} \varphi e^{2u_a}x_i\ud V_{g_{\S^2}}.
\end{align*}
Here we have used the fact $u_a \in \mathcal{S}$. 

To simplify
\begin{align*}
	\lambda_i:=3\fint_{\S^2} \left(-a\Delta_{\S^2} u_a+1\right)e^{-2u_a}x_i\ud V_{g_{\S^2}},
\end{align*}
we obtain
\begin{align*}
-a \Delta_{\S^2} u_a+1=e^{2u_a}\left(1+\sum_{i=1}^{3}\lambda_ix_i\right).
\end{align*}
Using the same argument in Chang and Yang \cite{Chang-Yang},  we conclude that $\lambda_i=0$ for $1\leq i \leq 3$. For a smaller $\ve_0$ if necessary, this together with Theorem \ref{Perturbation thm} implies $u_a= 0$, and also
\begin{align*}
\log \fint_{\S^2}  e^{2u} \ud V_{g_{\S^2}}\leq a\fint_{\S^2} |\nabla u|^2  \ud V_{g_{\S^2}}+2\fint_{\S^2} u  \ud V_{g_{\S^2}}
\end{align*}
for $a\in \left(\frac{1}{2}-\ve_0,\frac{1}{2}\right)$.
\hfill $\Box$

\vskip 16pt

\noindent \textbf{Data availability} Data sharing not applicable to this article as no datasets were generated or analysed during
the current study.

		\bibliographystyle{unsrt}

	\end{document}